\documentclass[11pt]{amsart}
\usepackage{amsfonts,amsmath,amssymb}
\usepackage[utf8]{inputenc}
\usepackage{fontenc}
\usepackage[all]{xy}
\usepackage{hyperref}

\begin{document}

\newcommand{\ba}{{\bf a}}
\newcommand{\bb}{{\bf b}}
\newcommand{\ab}{{\rm ab}}
\newcommand{\bpi}{{\boldsymbol{\pi}}}
\newcommand{\be}{{\bf e}}
\newcommand{\oh}{{\mathfrak o}}
\newcommand{\m}{{\mathfrak m}}
\newcommand{\jnf}{{\rm inf}}
\newcommand{\A}{{\mathbb A}}
\newcommand{\C}{{\mathbb C}}
\newcommand{\F}{{\mathbb F}}
\newcommand{\G}{{\mathbb G}}
\newcommand{\M}{{\mathbb M}}
\newcommand{\N}{{\mathbb N}}
\newcommand{\U}{{\mathbb U}}
\newcommand{\bP}{{\mathbb P}}
\newcommand{\PGl}{{\rm PGl}}
\newcommand{\R}{{\mathbb R}}
\newcommand{\bS}{{\mathbb S}}
\newcommand{\eL}{{\bf L}}
\newcommand{\el}{{\bf l}}
\newcommand{\Gal}{{\rm Gal}}
\newcommand{\Q}{{\mathbb Q}}
\newcommand{\bQ}{{\overline{\mathbb Q}}}
\newcommand{\T}{{\mathbb T}}
\newcommand{\SU}{{\rm SU}}
\newcommand{\Sp}{{\rm Sp \;}}
\newcommand{\W}{{\sf W}}
\newcommand{\sK}{{\sf K}}
\newcommand{\Z}{{\mathbb Z}} 
\newcommand{\Be}{{\boldsymbol{\epsilon}}}
\newcommand{\sslash}{{/\!/}}
\newcommand{\BHV}{{\rm BHV}}
\newcommand{\Spec}{{\rm Spec \; }}
\newcommand{\K}{{\rm K}}
\newcommand{\Rep}{{\rm R}}
\newcommand{\ox}{{\rm or}}
\newcommand{\bM}{{\overline{\mathbb M}}}
\newcommand{\B}{{\mathbb B}}
\newcommand{\cT}{{\match T}} 
\newcommand{\Maps}{{\rm Maps}}
\newcommand{\D}{{\mathbb D}}
\newcommand{\Con}{{\rm Con}}
\newcommand{\Diff}{{\rm Diff}}
\newcommand{\mM}{{\mathbb M}}
\newcommand{\Met}{{\rm Metrics}}
\newcommand{\mero}{{\rm mero}}
\newcommand{\Fns}{{\rm Fns}}
\newcommand{\Fred}{{\rm Fred}}
\newcommand{\Hom}{{\rm Hom}}
\newcommand{\Pic}{{\rm Pic}}
\newcommand{\wzw}{{\bf wzw}}
\newcommand{\FIO}{{\rm FIO}}
\newcommand{\boldelta}{{\boldsymbol{\Delta}}}
\newcommand{\Fin}{{\rm Fin}}
\newcommand{\cA}{{\mathcal A}}

\parindent=0pt
\parskip=6pt

\newcommand{\ie}{\textit{ie}\,}
\newcommand{\eg}{\textit{eg}\,}
\newcommand{\se}{{\sf e}}
\newcommand{\cf}{{\textit{cf}\,}}

\title{ Some very low-dimensional algebraic topology}

\author[J Morava]{J Morava}

\address{Department of Mathematics, The Johns Hopkins University,
Baltimore, Maryland} 

\email{jmorava1@jhu.edu}

\begin{abstract}{The Euclidean renormalization bundle considered in QFT by Connes, Kreimer, and Marcolli \cite{5,6} has been extended, in a remarkable series of papers by S Agarwala, to Riemannian manifolds $(X,g)$: in particular by the construction of a flat connection on that bundle \cite{1}(Prop 4.4), regarded as defined over a thickening $\Delta^\times \times X$ by an infinitesimal disk.}

\medskip \noindent The theory of Fourier integral operators on manifolds \cite{8,18} reconciles dimensional and zeta-function \cite{9} regularization by interpreting the disk $\Delta^\times$ as the germ of a neighborhood of a Jordan curve around $\infty \in \C P^1$. Such fields $X \to \Omega S^2$ were proposed in \cite{14} as useful in these contexts.\end{abstract}
  
\maketitle

{\bf \S I  Fields of finite subsets of the plane} \bigskip

{\bf 1.1} Following D Rolfsen and myself \cite{15}, Smale's (contractible) group $\D$ of compactly supported diffeomorphisms of the plane acts diagonally on the configuration  space 
\[
(\C^k - \boldelta^k)/{\mathfrak S}_k := \Con^k (\C)
\]
of unordered $k$-tuples of points in the plane. The graded topological transformation category
\[
\coprod_{k \geq 0} [\Con^k(\C) \sslash \D] := \Fin^\C_*
\]
(of finite subsets of the plane) has $\pi_0 = \N$ and the braid monoid $\B_*$ for its fundamental groupoid. \bigskip

With respect to Boardman-Vogt / Segal remote/disjoint union \cite{15}(\S 2.2 Prop 2) this is an associative monoid, with a translation operator (adding a remote point) which maps to a central element in the monoid ring $H_0 |\Fin^\C_*|$. 

After inverting this element (\ie stabilization), the group completion \cite{15}(\S 6)  
\[
|\Fin^\C_\otimes|^+ := \Omega B {\;} \Fin^\C_\otimes  \simeq \Omega S^2  \simeq B (\Omega^2 S^2)
\]
is then a moduli space for principal $\Omega^2 S^2$ - bundles, which we propose to call $\el$ines (or fields of finite subsets of the Dirac Sea \cite{14}(appendix {\it i})). Such things lie to the left of complex line bundles
\[
[\Omega S^2 \to B\T] \in H^2(\Omega S^2,\Z) = \Pic_\C(\Omega S^2)
\]
via the characteristic class Kronecker dual to $x^2 \in H_2(\Omega S^2,\Z)$. 

This group-like space is not simply-connected; indeed
\[
\pi_1 |\Fin^\C_\otimes|^+ \cong \Z 
\]
\cite{15}(Prop 9), and because loop spaces preserve fibrations, its universal cover
\[
\widetilde{|\Fin^\C_\otimes|^+} \simeq \Omega \SU(2) \sim (\Omega \Sigma) (S^2) \simeq \vee_{k \geq 0} S^{2k} \simeq
S^0[S^2]
\]
looks stably like a polynomial algebra on one generator.\bigskip

Similarly, $\Omega^2 \SU(2) \to \Omega^2_e \C P^1$ is a homotopy equivalence, so we might think of a $\el$ine (bundle) over $X$ as an $\Omega^2 \R^3_+$-torsor or field of Wess-Zumino-Witten models. In the terminology of \cite{14}(appendix {\it xiii} \S 1.3), the immersions in $\Omega S^2$ define a space of classical states of the Nambu-Goto model : perturbations of the equator $S^1 \subset S^2$, smooth Jordan curves encircling $0 \in \C P^1$, wriggles of the Midgard serpent \dots \bigskip

{\bf 1.2 example}  The ring homomorphism
\[
{\mathbb H} \to {\mathbb H} \otimes_\R \C \cong  {\mathbb M}_2(\C)
\]
defines a group homomorphism 
\[
[\SU(2) \to \Maps (\C P^1, \C P^1) \subset \Omega^2_e S^2 ] \in \pi_3 \Omega^2 S^2
\]
on unit groups analogous to the Bott-Cartan class in a simple Lie group \eg $\pi_3 \SU(2)$. Its delooping 
\[
\eta^2  = [S^1 \wedge S^3  \to  S^2] \in \pi_4 S^2 = \Z_2
\]
is the first element of coker $J$: recall \cite{16}(\S  1.1.12) that $\pi^S_1 \ni \eta \mapsto \{\pm\} \in K_1\Z $, and that  $\eta^3 \equiv 4\nu$ mod 2. This defines a canonical class

\[
\xymatrix{
B(\SU(2) \times \SU(2)) \cong B{\rm Spin}(4) \ar[d]^{B(\eta^2 \times \eta^2)} \\
B(\Omega^2 S^2) \simeq \Omega S^2 }
\]
in 
\[
\Pic_\wzw B{\rm Spin}(4) = \pi_0 \Maps (B{\rm Spin}(4),B\Omega^2_e S^2) ,
\]

perhaps of interest \cite{13}(\S 2.5) in Hopkins and Singer's study \cite{10}(appendix E) of relative deformations of spin structures \dots 

\newpage

{\bf 1.3} A theorem of Caratheodory \cite{20} identifies elements
\[
\{ b \in H^\infty(\C_{<1}) {\; | \:} b \subset \bb : \C_{\leq 1} \to \C\} 
\]
of the Hardy space (of holomorphic functions defined and bounded on $\C_{<1}$) which are conformal isomorphisms on the interior and extend to the boundary, with disks in the plane bounded by Jordan curves encircling the origin. This defines a ($\T$-equivariant) map $b \mapsto \partial b := \bb|_{\C_{=1}}$
\[
\xymatrix{ 
{} & \Omega S^2 \ar[d]  \\
\{\rm Jordan\} \ar@{.>}[ur] \ar[r]^\partial  & L\C_+}
\]
from the space of such disks to the space of free loops on $S^2$. A contour around the boundary circle and down the image of $t \in i\R_{\geq 0}$ defines a loop based at 0, lifting to a map to $\Omega S^2$.  Sewing along the outgoing boundary\begin{footnote}{$\exp(it) = z \mapsto (z^{-1}  + \sum_{k \geq 2} b_k z^{-k})^{-1} \equiv z {\; \rm mod \;} \C[[z^{-1}]]$ \cite{19} Prop 2.2}\end{footnote} defines an action 
\[
\cA \times \{\rm Jordan\} \to \{\rm Jordan\}
\]
of Segal's monoid $\cA$ of annuli, making \{\rm Jordan\} an $\cA$-line (roughly, free of rank one), suggesting that renormalizability entails an underlying homotopy-theoretical proto-CFT. \bigskip

{\bf \S II {\;} $\phi^4$ GR/QED} following S Agarwala \bigskip

{\bf 2.1} On a Riemannian background $(X,g)$, BCKM renormalization \cite{5, 6, 14} associates to a section (equivalently, a trivialization)

\[
\xymatrix{
\C P^1 \times X \ar[r]^{\rm Seeley} & P \ar[d]^{\langle \rm Feynman-Lie \rangle} \ar[r]^{{\rm eval}_{U(\Gamma)}} & \C \\
{} & \ar[ul] \Delta^\times \times X \ar@{.>}[ur] }
\] 

of a certain renormalization bundle (with fiber a graded prounipotent group associated to a combinatorial Lie algebra of Feynman diagrams \cite{11} under insertion), together with a point (valued in something like $z^{-1}\C[[z]]$) of that bundle, a renormalized value \cite{6}(Th 2.5), \cite{2}(\S 4.3) in $\C$ for a Feynman-Wick integral $U(\Gamma)$ in a physical model, \cite{7}, \cite{6}(\S 2.7). Agarwala's flat connection \cite{1}(Props 4.4, 4.6) on the renormalization bundle implies that two trivializations yield gauge equivalent renormalized values. 

The adjoint composition
\[
X \to \Maps_\mero(\C P^1,P) \to \Maps (\Delta^\times,P) \to \C
\]

factors through an identification of the infinitesimal disk $\Delta^\times$ with a neighborhood of $\infty \in \C P^1$ \cite{18}(\S 2.2), reinterpreting \cite{3}(Th 3.9) the dimensional regularization parameter in the Riemannian context in terms of the complex power formalism of \cite{8,9,18}. 

In this picture, the renormalized value of the integral $U(\Gamma)$ indexed by a diagram $\Gamma$ is the Cauchy residue associated to a Jordan curve defined by a field $X \to \Omega S^2$ \dots \bigskip

{\bf 2.2} Suppose then that $(X,\partial X = Y)$ is a compact connected oriented smooth four-dimensional manifold; if $\pi_0 Y$ is nontrivial, fix open book decompositions $\Sigma_i \to Y_i - L_i \to S^1$ (\ie bound by links $L_i \subset Y_i$) for its components.

Let $\D := \Diff (X |  \cdots)$ be the group(oid) of diffeomorphisms  of $X$ (subject to conditions $ | \cdots$), and similarly let $\mM := \Met (X | \cdots)$ be the (contractible) space of Riemanian metrics on $X$ (\eg with boundary conditions) such that $\D$ acts cromulently on $\M$, \cf [Ebin]. The homotopy quotient 
\[
\mM_{h\D} = \mM \times_\D  E \D \to {\rm pt} \times_\D E \D = B \Diff (X|\cdots)
\]

has compact Lie groups as isotropy, making Spec $H_\D^*(\mM)$ accessible \dots \bigskip

{\bf 2.3}  An idealized model for a renormalization scheme following Marcolli for classical QED \cite{17} and Agarwala for its general relativisation might be an $\infty$-functor
\[
[\C_+ \times \Met \sslash \D ] \Rightarrow [\Fns_\mero (X,\FIO) \sslash \D]
\]
exemplified by a system of Seeley (meromorphic) traces
\[
(s, g) \mapsto  {\rm Tr} (-\Delta_g)^{-s}
\]
or more generally by Fourier integral operators $U(\Gamma)$ \cite{6}(\S 2.2.2, Th 2.5-6, \S 2.7), \cite{8}(VI \S 6) associated to configuration spaces in $\R^4$ indexed by Feynman graphs. \bigskip

Quasi-shuffle zeta-function renormalization following [Cartier, Hoffman, \cite{14}(\S 2.3)] provides a basic formalism for a Haag-Kastler approach to QED/GR on 4-dimensional Riemannian cobordisms between three-manifolds. \bigskip

[in progress: some further explanatory material is included in the .tex file]

\newpage 

\bibliographystyle{amsplain}

\end{document}